\documentclass[9pt, final, oneside, twocolumn, letter]{autart_modified}

\usepackage[authoryear]{natbib}
\usepackage[draft]{hyperref}

\usepackage{graphicx}
\usepackage{amsmath,mathtools} %EPS Graphics
\usepackage{amssymb, amsmath}
\usepackage{color}
\usepackage{caption}
\usepackage{tikz}
\usetikzlibrary{positioning}
\usetikzlibrary{shapes,arrows}
\usepackage{epstopdf}

\usepackage{enumitem,kantlipsum}
\usepackage{mathtools,lipsum,cuted}
\newcommand{\lls}{\langle\langle}
\newcommand{\ggs}{\rangle\rangle}

\newtheorem{definition}{Definition}
\newtheorem{corollary}{Corollary}

\newtheorem{remark}{Remark}
\newtheorem{proposition}{Proposition}

\date{}

\begin{document}
\begin{frontmatter}

\title{A Central Synergistic Hybrid Approach for Global Exponential Stabilization on $SO(3)$}
\thanks[footnoteinfo]{This work was supported by the National Sciences and Engineering Research Council of Canada (NSERC).}

\author[First]{Soulaimane Berkane}\ead{sberkane@uwo.ca}
\author[First]{Abdelkader Abdessameud}\ead{aabdess@uwo.ca}
\author[First,Second]{Abdelhamid Tayebi}\ead{atayebi@lakeheadu.ca}

\address[First]{Department of Electrical and Computer Engineering, University of Western Ontario, London, Ontario, Canada, N6A 3K7.}
\address[Second]{Department of Electrical Engineering, Lakehead University, Thunder Bay, Ontario, Canada, P7B 5E1.}

\begin{keyword}
Attitude control; Hybrid feedback; Synergistic potential functions; Special orthogonal group.
\end{keyword}

\begin{abstract}
We propose a new central synergistic hybrid approach for global exponential stabilization on the Special Orthogonal group $SO(3)$. We introduce a new switching concept relying on a central family of (possibly) non-differentiable potential functions that enjoy (as well as their gradients) the following properties: 1) being quadratic with respect to the Euclidean attitude distance, and 2) being synergistic with respect to the gradient's singular and/or critical points. The proposed approach is used to solve the attitude tracking problem, leading to global exponential stability results.
\vspace{-0.2 in}
\end{abstract}
\end{frontmatter}
\section{Introduction}
            The rigid body attitude control problem has received a growing interest during the last decade, with various applications in aerospace, marine engineering, and robotics \citep[see, for instance,][to name a few]{hughes1986sad, joshi1995robust, tayebi2006attitude,seo2007separation,sakagami2010attitude, kreutz1988attitude}.
            Attitude control schemes can be categorized by the choice of the attitude parameterization, such as Euler-angles, unit quaternion, and modified Rodrigues parameters. All these parameterizations fail to represent the attitude of a rigid body both globally and uniquely, which results in control schemes that are either singular or exhibit the so-called unwinding phenomenon \citep{Chaturvedi2011}. On the other hand, some research work such as \citep{Koditschek,suzuki2005observer,Sanyal2009, lee2012,bayadi2014almost}, has been devoted to the attitude control design directly on the Special Orthogonal group $SO(3)$, where the orientation of the rigid body is uniquely and globally represented by a rotation matrix in $SO(3)$.

            %such as in \cite{Koditschek, Sanyal2009, lee2012}, has been devoted to the attitude control and estimation design directly on the Special Orthogonal group $SO(3)$, where the orientation of the rigid body is uniquely and globally represented by a rotation matrix in $SO(3)$.
            Indeed, the group $SO(3)$ has the distinct feature of being a boundaryless compact manifold with a Lie group structure that allows the design and analysis of attitude control systems within the well established framework of geometric control \citep{bullo2005geometric}.
\vspace{0.05 in}\\
             The group $SO(3)$ is not diffeomorphic to any Euclidean space and hence there does not exist any continuous time-invariant feedback on $SO(3)$ that achieves global asymptotic stability \citep{Bhat2000}. In \cite{Koditschek}, for instance, a continuous time-invariant control scheme has been shown to asymptotically track any smooth reference attitude trajectory starting from arbitrary initial conditions except from a set of Lebesgue measure zero. This is referred to as \textit{almost} global asymptotic stability, and is mainly due to the appearance of undesired critical points (equilibria) when using the gradient of a smooth potential function in the feedback. In fact, any smooth potential function on $SO(3)$ is guaranteed to have at least four critical points where its gradient vanishes \citep{morse1934calculus}.
\vspace{0.05 in}\\
           In \cite{mayhew2011hybrid}, a hybrid feedback scheme has been proposed to overcome the topological obstruction to global asymptotic stability on $SO(3)$ and, at the same time, ensure some robustness to measurement noise. The main idea in the latter paper is to design a hybrid algorithm based on a family of \textit{smooth} potential functions and a hysteresis-based switching mechanism that selects the appropriate control action corresponding to the minimal potential function.
%to avoid the undesired critical points. After each switching, the control law derived from the minimal potential function is selected. A sufficient condition for global asymptotic stability of the resulting hybrid controller is the ``synergism" property.
           It was shown that a sufficient condition to avoid the undesired critical points, and ensure global asymptotic stability, is the ``synergism" property (given in Definition~\ref{definition::synergism} later) of the smooth potential functions. The angular warping technique has been used in \cite{mayhew2011synergistic}, and extended later in \cite{casau2014globally} and \cite{berkane2015construction}, to construct \textit{central} synergistic potential functions on $SO(3)$.
           %Note that, in \cite{Mayhew2013} and \cite{lee2015tracking}, different approaches have been adopted to construct \textit{non-centrally} synergistic potential functions on $SO(3)$ from which globally asymptotically stabilizing control laws can be designed.
           In \cite{Mayhew2013} and \cite{lee2015tracking}, different approaches have been adopted to construct a \textit{non-central} synergistic family of smooth potential functions from which globally stabilizing control laws can be designed. It is important to mention here that, in contrast to the non-central approach, the control algorithm derived from each potential function in the central synergistic family guarantees (independently) almost global asymptotic stability results. It is also worth pointing out that non-central and central synergistic potential functions have been considered in \cite{mayhew2013spherical} and \cite{casau2015global} to ensure, respectively, global asymptotic and global exponential stabilization on the $n$-dimensional sphere. However, the extension of these approaches to the \textit{full} attitude control problem on $SO(3)$ is not straightforward.
\vspace{0.05 in}\\
           In this paper, we propose solutions to the global exponential stabilization problem on $SO(3)$ using hybrid feedback. Specifically, our contributions can be summarized as follows.
\vspace{-0.1 in}
        \begin{itemize}[wide, labelwidth=!, labelindent=0pt]
        \item  We introduce a new \textit{exp-synergism} concept for hybrid attitude control systems design on $SO(3)$. This extends the central synergism concept introduced in \cite{mayhew2011synergistic, Mayhew2013} by imposing more restrictive conditions on the class of potential functions (not necessarily differentiable everywhere on $SO(3)$) used in the attitude control design.
        As such, significant advantages are gained such as exploiting non-smooth potential functions, which have been shown to ensure better performance compared to the traditional smooth potential functions \cite{lee2012,zlotnik2015nonlinear}.
       % We extend the central synergism concept introduced in \cite{mayhew2011synergistic,Mayhew2013} by imposing more restrictive conditions on the potential functions used in the attitude control design. The new concept, named \textit{exp-synergism}, allows the use of a wider class of potential functions that are not necessarily differentiable everywhere on $SO(3)$. This provides significant benefits when exploiting non-smooth potential functions, which have been shown to provide faster convergence rates for large attitude errors compared to the traditional smooth potential functions \cite{lee2012,zlotnik2015nonlinear}.
%and define a new concept, namely \textit{exp-synergism}, which imposes more restrictive conditions on the potential function used in the attitude control design and analysis. In contrast to the classical synergism concept, our concept of exp-synergism allows the use of non-differentiable and/or differentiable potential functions on $SO(3)$. This allows to take advantage of the large class of potential functions used in the literature, which may not be differentiable everywhere on $SO(3)$. It is well known that non-smooth potential functions are likely to provide a faster convergence rate for large attitude errors compared to the traditional trace functions used often in the literature \cite{lee2012,zlotnik2015nonlinear}.
         \item We propose hybrid control schemes guaranteeing global exponential stabilization of the single and double integrator systems on $SO(3)$ and $SO(3)\times\mathbb{R}^3$, respectively. The proposed schemes employ a switching mechanism that switches between different configurations of exp-synergistic potential functions. In contrast to \cite{mayhew2011hybrid}, the hybrid switching allows to avoid a singular 3-dimensional manifold (instead of a finite number of critical points) where the gradient of the potential function is not defined. In contrast to \cite{lee2015tracking, Mayhew2013}, the proposed family of exp-synergistic potential functions is central which, as mentioned above, guarantees that each configuration of the hybrid closed-loop system ensures  almost global asymptotic stability.
         \item For each of the proposed control schemes, we provide an alternative design ensuring a smooth control input without sacrificing the global exponential stability result. The proposed control \textit{smoothing} procedure is much simpler than the backstepping procedure proposed in \cite{Mayhew2013}.
         \item As an application of our approach, we present a solution to the rigid body attitude tracking control problem guaranteeing global exponential stability. To the best of our knowledge, there is no similar solution to the full attitude tracking control problem in the literature.
\end{itemize}
\section{Preliminaries and Problem Formulation}\label{sec2}
\vspace{-0.05 in}
\subsection{Notations}\label{notation}
\vspace{-0.1 in}
          Throughout the paper, we use $\mathbb{R}$, $\mathbb{R}^+$ and $\mathbb{N}$ to denote, respectively, the sets of real, nonnegative real and natural numbers. We denote by $\mathbb{R}^n$ the $n$-dimensional Euclidean space and by $\mathbb{S}^n$ the unit $n$-sphere embedded in $\mathbb{R}^{n+1}$. Given $x\in\mathbb{R}^n$, the Euclidean norm of $x$ is defined as $\|x\|=\sqrt{x^\top x}$. For a square matrix $A\in\mathbb{R}^{n\times n}$, we denote by $\lambda_i^A, \lambda_{\mathrm{min}}^A$, and $\lambda_{\mathrm{max}}^A$ the $i$th, minimum, and maximum eigenvalue of $A$, respectively. Given two matrices $A, B\in\mathbb{R}^{m\times n}$, their Euclidean inner product is defined as $\lls A,B\ggs=\textrm{tr}(A^{\top}B)$ and the Frobenius norm of $A$ is $\|A\|_F=\sqrt{\lls A,A\ggs}$.  A function $f:X\to Y$ is a $\mathcal{C}^k$ function if all its first $k$ derivatives exist and are continuous, in which case we write $f\in\mathcal{C}^k(X,Y)$.

\vspace{-0.05 in}
\subsection{Attitude representation}
\vspace{-0.1 in}
            The orientation (or attitude) of a rigid body is represented by a {\it rotation matrix} $R$ belonging to the special orthogonal group $SO(3) := \{ R \in \mathbb{R}^{3\times 3}|\; \mathrm{det}(R)=1,\; RR^{\top}= I \}$, where $I$ is the three-dimensional identity matrix. The group $SO(3)$ has a compact manifold structure with its tangent spaces being identified by $T_RSO(3):=\left\{R\Omega\mid\Omega\in\mathfrak{so}(3)\right\}$. The \textit{Lie algebra} of $SO(3)$, denoted by $\mathfrak{so}(3):=\left\{\Omega\in\mathbb{R}^{3\times 3}\mid\;\Omega^{\top}=-\Omega\right\}$, is the vector space of 3-by-3 skew-symmetric matrices. The inner product on $\mathbb{R}^{3\times 3}$, when restricted to the Lie-algebra of $SO(3)$, defines the following \textit{left-invariant} Riemannian metric on $SO(3)$
            \begin{eqnarray}\label{metric}
            \langle R\Omega_1,R\Omega_2\rangle_R:=\lls\Omega_1,\Omega_2\ggs,
            \end{eqnarray}
            for all $R\in SO(3)$ and $\Omega_1,\Omega_2\in\mathfrak{so}(3)$. The map $[\cdot]_\times: \mathbb{R}^3\to\mathfrak{so}(3)$ is defined such that $[x]_\times y=x\times y$, for any $x, y\in\mathbb{R}^3$, where $\times$ denotes the vector cross-product on $\mathbb{R}^3$. Also, let $\mathrm{vex}:\mathfrak{so}(3)\to\mathbb{R}^3$ denote the inverse isomorphism of the map $[\cdot]_\times$, such that $\mathrm{vex}([\omega]_\times)=\omega,$ for all $\omega\in\mathbb{R}^3$ and $[\mathrm{vex}(\Omega)]_\times=\Omega,$ for all $\Omega\in\mathfrak{so}(3)$. Defining $\mathbb{P}_a:\mathbb{R}^{3\times 3}\to\mathfrak{so}(3)$ as the projection map on the Lie algebra $\mathfrak{so}(3)$ such that $ \mathbb{P}_a(A):=(A-A^{\top})/2$, one can extend the definition of $\mathrm{vex}$ to $\mathbb{R}^{3\times 3}$ by taking the composition map $\psi := \mathrm{vex}\circ \mathbb{P}_a$ %$:\mathbb{R}^{3\times 3}\to\mathbb{R}^3$
            such that, for a $3$-by-$3$ matrix $A:=[a_{ij}]_{i,j = 1,2,3}$, one has
\vspace{-0.1 in}
            \begin{eqnarray}\label{psi}
            \psi(A):=\mathrm{vex}\left(\mathbb{P}_a(A)\right)=\frac{1}{2}\setlength{\arraycolsep}{6pt}
                        \renewcommand{\arraystretch}{0.2}\begin{bmatrix}
            a_{32}-a_{23}\\a_{13}-a_{31}\\a_{21}-a_{12}
            \end{bmatrix}.\vspace{-0.05 in}
            \end{eqnarray}
            The attitude of a rigid body can also be represented as a rotation of angle $\theta\in\mathbb{R}$ around
            a unit vector axis $u\in\mathbb{S}^2$. This is known as the angle-axis parametrization of $SO(3)$ and is given by the map $\mathcal{R}:\mathbb{R}\times\mathbb{S}^2\to SO(3)$ such that
\vspace{-0.1 in}
            \begin{eqnarray}
            \mathcal{R}(\theta,u)=I+\sin(\theta)[u]_\times+(1-\cos\theta)[u]_\times^2.\nonumber
            \end{eqnarray}
            For any attitude matrix $R\in SO(3)$, we define $|R|_I\in[0, 1]$ as %the normalized Euclidean distance on $SO(3)$ which is given by %\vspace{-0.1 in}
            $$|R|_I^2:=\mathrm{tr}(I-R)/4.$$
%\vspace{-0.1 in}
\subsection{Exp-synergistic potential functions on SO(3)}\label{section::synergistic}
\vspace{-0.1 in}
            Given a finite index set $\mathcal{Q}\subset\mathbb{N}$, we let $\mathcal{C}^0\left(SO(3)\times\mathcal{Q},\mathbb{R}^+\right)$ denote the set of positive-valued and continuous functions $\mathcal{U}: SO(3)\times\mathcal{Q}\to\mathbb{R}^+$. If, for each $q\in\mathcal{Q}$, the map $R\mapsto\mathcal{U}(R,q)$ is differentiable on the set $D_q\subseteq SO(3)$ then the function $\mathcal{U}(R,q)$ is continuously differentiable on $\mathcal{D}\subseteq SO(3)\times\mathcal{Q}$, where
            $\mathcal{D}=\cup_{q\in\mathcal{Q}}D_q\times\{q\}$, in which case we denote $\mathcal{U}\in\mathcal{C}^1\left(\mathcal{D},\mathbb{R}^+\right)$. Additionally, for all $(R,q)\in\mathcal{D}$, $\nabla\mathcal{U}(R,q)\in T_RSO(3)$ denotes the gradient of $\mathcal{U}$, with respect to $R$, relative to the Riemannian metric \eqref{metric}.%\\
\vspace{0.05 in}\\
            A function $\mathcal{U}\in\mathcal{C}^0\left(SO(3)\times\mathcal{Q},\mathbb{R}^+\right)$ is said to be  a potential function on $\mathcal{D}\subseteq SO(3)\times\mathcal{Q}$ with respect to the set $\mathcal{A}=\{I\}\times\mathcal{Q}$ if: $i)$ $\mathcal{U}(R,q)>0$ for all $(R,q)\notin\mathcal{A}$, $ii)$ $\mathcal{U}(R,q)=0$, for all $(R,q)\in\mathcal{A}$, and $iii)$ $\mathcal{U}\in\mathcal{C}^1(\mathcal{D},\mathbb{R}^+)$.
%\vspace{-0.1 in}
%            \begin{itemize}
%            \item $\mathcal{U}(R,q)>0$ for all $(R,q)\notin\mathcal{A}$,
%            \item $\mathcal{U}(R,q)=0$, for all $(R,q)\in\mathcal{A}$,
%            \item $\mathcal{U}\in\mathcal{C}^1(\mathcal{D},\mathbb{R}^+)$.
%         \vspace{-0.1 in}
%             \end{itemize}
%
             The set of all potential functions on $\mathcal{D}$ with respect to $\mathcal{A}$ is denoted as $\mathcal{P}_\mathcal{D}$, where a function $\mathcal{U}(R,q)\in\mathcal{P}_\mathcal{D}$ can be seen as a family of potential functions on $SO(3)$ encoded into a single function indexed by $q$.
\vspace{-0.1 in}
        \begin{definition}[Synergism]\cite{Mayhew2013}\label{definition::synergism}
            For a given finite index set $\mathcal{Q}\subset\mathbb{N}$, let $\mathcal{U}\in\mathcal{P}_{SO(3)\times\mathcal{Q}}$. The potential function $\mathcal{U}$ is said to be centrally synergistic, with gap exceeding $\delta>0$, if and only if the condition \vspace{-0.05 in}
                \begin{eqnarray}\label{condition::synergistic}
				&&\mathcal{U}( R,q)-\underset{m\in\mathcal{Q}}{\mathrm{min}}\;\mathcal{U}(R,m)>\delta, \vspace{-0.2 in}
               \end{eqnarray}
               holds for all $(R,q)\in SO(3)\times\mathcal{Q}$ such that $\nabla\mathcal{U}(R,q)=0$ and $(R,q)\notin\mathcal{A} = \{I\}\times \mathcal{Q}$.
\end{definition}
            Condition \eqref{condition::synergistic} implies that at any given critical point $(R,q)$, there exists another point $(R,m)\in SO(3)\times\mathcal{Q}$ such that $\mathcal{U}(R,m)$ has a lower value than $\mathcal{U}(R,q)$. The property in Definition~\ref{definition::synergism} has allowed to design attitude control systems with global asymptotic stability results \cite{mayhew2011hybrid,Mayhew2013}. The term ``central" refers to the fact that all the potential functions $R\mapsto\mathcal{U}(R,q)$ share the identity element $I$ as a critical point such that $\nabla\mathcal{U}(I, q)=0$ for all $q\in\mathcal{Q}$. %, \textit{i.e.} one has $\nabla\mathcal{U}(\mathcal{A})=0$.
            %Non-central synergistic potential functions have also been used in \cite{Mayhew2013} and \cite{lee2015tracking} as mentioned in the Introduction section.
%to achieve global asymptotic stability (GAS) and global exponential stability (GES). However, in contrast to the non-central approach, each controller derived from a potential function in a central synergistic family guarantees (independently) almost global asymptotic stability.
        \begin{definition}[Exp-synergism]\label{definition::expsynergism}
            For a given finite index set $\mathcal{Q}\subset\mathbb{N}$ and $\mathcal{D}\subseteq SO(3)\times\mathcal{Q}$, let $\mathcal{U}\in\mathcal{P}_\mathcal{D}$ and the set
            $C\subseteq SO(3)\times\mathcal{Q}$ defined as \vspace{-0.05 in}
                \begin{eqnarray}
               C:=\{(R,q): \mathcal{U}( R,q)-\underset{m\in\mathcal{Q}}{\mathrm{min}}\;\mathcal{U}(R,m)\leq\delta\}. \vspace{-0.1 in}
                \end{eqnarray}
            The potential function $\mathcal{U}$ is said to be
                {\it exp}-synergistic, with gap exceeding $\delta>0$, if and only if there exist constants $\alpha_i>0$, $i=1,\ldots, 4$, such that the following hold:
                \begin{eqnarray}\label{condition::expsynergistic1}
				&&\alpha_1|R|_I^2\leq\mathcal{U}(R,q)\leq\alpha_2|R|_I^2,\quad (R,q)\in SO(3)\times\mathcal{Q},\\
                \label{condition::expsynergistic2}
                &&\alpha_3|R|_I^2\leq\left\|\nabla\mathcal{U}(R,q)\right\|^2_F\leq\alpha_4|R|_I^2, \quad (R,q)\in C,\\\label{condition::expsynergistic3}
               &&C\subseteq\mathcal{D}.
               \end{eqnarray}
\end{definition}
         \vspace{-0.15 in}
                Definition \ref{definition::expsynergism} considers a wider class of potential functions $\mathcal{U}\in\mathcal{P}_\mathcal{D}$ and imposes more restrictive conditions as compared to Definition \ref{definition::synergism}. The exp-synergism property, as will become clear later, plays an important role to ensure desirable exponential decay when using a gradient-based feedback on $SO(3)$.  It can be verified that if $\mathcal{U}\in\mathcal{P}_{SO(3)\times\mathcal{Q}}$ is an exp-synergistic potential function then it is synergistic as well. In fact, condition \eqref{condition::expsynergistic2} implies that the gradient $\nabla\mathcal{U}(R,q)$ does not vanish expect at $\mathcal{A}$ which is equivalent to condition \eqref{condition::synergistic}. The opposite does not hold in general.
\vspace{-0.05 in}
        \subsection{Models and problem statement}\label{Sec:models}
\vspace{-0.1 in}
            Consider the following system on $SO(3)$ \vspace{-0.1 in}
                \begin{eqnarray}\label{dR:1}
                    \dot R=R[u_1]_\times, \vspace{-0.1 in}
                \end{eqnarray}
            where $R\in SO(3)$ represents the attitude state and $u_1\in\mathbb{R}^3$ is the control input. The stabilization control problem for system \eqref{dR:1} consists in designing appropriate input $u_1$ such that $(R=I)$ is globally exponentially stable.
%
%\vspace{0.05cm}\\
%
            We also consider system \eqref{dR:1} augmented by the angular velocity dynamics
                \begin{eqnarray}\label{dR:2}
                \dot R=R[\omega]_\times,\qquad \dot\omega = u_2, \vspace{-0.1 in}
                \end{eqnarray}
            where $R\in SO(3)$ and $\omega\in\mathbb{R}^3$ represent the state variables and $u_2\in\mathbb{R}^3$ is the control input. The objective is to design an appropriate input $u_2$ such that $(R=I, \omega=0)$ is globally exponentially stable.
%
%\vspace{0.05cm}\\
%
%For practical purposes, we also propose a smoothing mechanism that removes the switching effects from the actual control input through an integrator.
 \vspace{-0.1 in}
\section{Main results}\label{Sec:mainresult}
\vspace{-0.1 in}
            In this section, we propose global exponential hybrid attitude control schemes on $SO(3)$ for both systems \eqref{dR:1} and \eqref{dR:2}. For the sake of clarity, we present our main results in two parts.
            First, we propose an approach for the design of hybrid control algorithms using a switching dynamical system based on the gradient of generic potential functions on $SO(3)$. Basically, we
%show that global exponential stability of the closed loop equilibrium point can be guaranteed with a particular choice of potential functions satisfying some sufficient conditions.
            derive the sufficient conditions on these potential functions such that global exponential stability of the closed-loop equilibrium point is guaranteed. We also show that the same results are guaranteed using continuous-time control laws that are obtained via a dynamic extension of the originally derived hybrid controllers.
            %to smoothen the can be  the attitude proportional term (which contains discontinuous jumps) into a first order low-pass filter allows to obtain a continuous control input guaranteeing the same exponential results.
             Second, we provide a systematic methodology for the construction of such potential functions on $SO(3)$, and derive the expression of the control algorithms accordingly.

\vspace{-0.05 in}
\subsection{Global Exponential Feedback Control on $SO(3)$}\label{Sec:sub:exponential}
\vspace{-0.1 in}
        In this subsection, we show that exp-synergistic potential functions are instrumental for the global exponential stabilization of \eqref{dR:1} and \eqref{dR:2}.
%Moreover, these hybrid controllers can be made continuous by adding a low pass filter without sacrificing the GES results.
%
%Inspired by \cite{mayhew2011hybrid}, the control input is based on a hybrid mechanism that uses a ``min-switch" strategy to select the control law derived from the minimal potential function among some family of potential functions on $SO(3)$.
        Let $\mathcal{Q}\subset \mathbb{N}$ be a finite index set and $\mathcal{U}\in\mathcal{P}_\mathcal{D}$ be a potential function on $\mathcal{D}\subseteq SO(3)\times\mathcal{Q}$. For a given $\delta>0$, the switching mechanism of the discrete state variable $q$, which dictates the current \textit{mode} of operation of the hybrid control system\footnote{In this work, we make use of the recent framework for dynamical hybrid systems found in \citet{goebel2012hybrid}.}, is given by
                \begin{eqnarray}\label{q}
                \left\{
                \begin{array}{ll}
                \dot q=0,\;\;\;&(R,q)\in C,\\
                q^+\in \mathrm{arg}\min_{m\in\mathcal{Q}}\;\mathcal{U}(R,m),\;\;\;&(R,q)\in D,
                \end{array}\right.
                \end{eqnarray}
        where the flow set $C$ and jump set $D$ are defined by
                \begin{eqnarray}
                \label{C}C&=&\{(R,q): \mathcal{U}( R,q)-\underset{m\in\mathcal{Q}}{\mathrm{min}}\;\mathcal{U}(R,m)\leq\delta\},\\
                \label{D}D&=&\{(R,q):  \mathcal{U}(R,q)-\underset{m\in\mathcal{Q}}{\mathrm{min}}\;\mathcal{U}(R,m)\geq\delta\}.
                \end{eqnarray}
        It is worthwhile mentioning that the above hybrid mechanism, inspired by \cite{mayhew2011hybrid}, uses a ``min-switch" strategy to select a control law derived from the minimal potential function among some family of potential functions on $SO(3)$. Such a control input for systems \eqref{dR:1} and \eqref{dR:2} is given, respectively, in Theorem~\ref{theorem::kinematics} and Theorem~\ref{theorem::dynamics} proved in the Appendix.
        %
        %Our results in this subsection are given in Theorem~\ref{theorem::kinematics}, Theorem~\ref{theorem::dynamics} and Theorem~\ref{theorem::dynamics::smooth} proved in Appendix \ref{proof::theorem::kinematics}, Appendix \ref{proof::theorem::dynamics} and Appendix \ref{proof::theorem::dynamics::smooth}, respectively.
%
   %
            \begin{thm}\label{theorem::kinematics}%(Hybrid P-type control)\\
            Consider system \eqref{dR:1} with the control input
            \begin{eqnarray}\label{control1}
                        u_1=-k_cx_R(R,q),\quad k_c>0,
            \end{eqnarray}
            with $x_R(R,q):=\psi\left(R^{\top}\nabla\mathcal{U}(R,q)\right)$, where $q$ is generated by \eqref{q}-\eqref{D} and  the map $\psi(\cdot)$ is defined in \eqref{psi}. If the potential function $\mathcal{U}\in\mathcal{P}_\mathcal{D}$ used in \eqref{q}-\eqref{control1} is {\it exp}-synergistic with synergistic gap exceeding $\delta$, then the set $\mathcal{A}=\{I\}\times\mathcal{Q}$ is globally exponentially stable.
            \end{thm}
            \begin{thm}\label{theorem::dynamics}%(Hybrid PD-type control)\\
            Consider system \eqref{dR:2} with the control input
             \begin{eqnarray}\label{control2}
             u_2=-k_cx_R(R,q)-k_\omega\omega,\quad k_c,~k_\omega>0,
            \end{eqnarray}
            where $x_R(R,q)$ is defined as in Theorem~\ref{theorem::kinematics}. If the potential function $\mathcal{U}\in\mathcal{P}_\mathcal{D}$ used in \eqref{q}-\eqref{D} and \eqref{control2} is {\it exp}-synergistic with synergistic gap exceeding $\delta$, then the set $\bar{\mathcal{A}}:=\{(R, q, \omega)\in SO(3)\times\mathcal{Q}\times\mathbb{R}^3\mid\;R=I, \omega = 0\}$ is globally exponentially stable.
            \end{thm}
            Note that the proposed attitude control schemes are subject to discontinuous jumps due to the direct switching in the hybrid controller configuration. However, the result of Theorem \ref{theorem::dynamics} can be used to design a continuous control input for \eqref{dR:1}, which is generally more desirable in practice. In fact, letting the input $u_1$ in \eqref{dR:1} be given by the output of the dynamic system
            \begin{eqnarray}\label{du1}
            \dot u_1=-k_\omega u_1-k_cx_R(R,q),\quad u_1(0)=u_{1,0}, %\in\mathbb{R}^3,
            \end{eqnarray}
            guarantees, by Theorem \ref{theorem::dynamics}, global exponential stability of the set $\mathcal{A}_{s}=\{(R,q,u_1)\in SO(3)\times\mathcal{Q}\times\mathbb{R}^3\mid\;R=I, u_1=0\}$. Using this observation, the attitude control scheme in Theorem \ref{theorem::dynamics} can be modified, by introducing a simple dynamic system, to ensure the same global exponential stability result using continuous control for system \eqref{dR:2}. This is shown in the following result proved in Appendix~\ref{proof::theorem::dynamics::smooth}.
\vspace{-0.1 in}
            \begin{thm}\label{theorem::dynamics::smooth}%(Continuous PD-type control)\\
            Consider system \eqref{dR:2} with the control input
            \begin{eqnarray}\label{control3}\begin{array}{lcl}
            u_2&=&-k_cx_{R,s}-k_\omega\omega,\\
            \dot x_{R,s}&=&-k_s(x_{R,s}-x_{R}(R,q)),\end{array}
            \end{eqnarray}
            for some $x_{R,s}(0)$, where $k_c,k_\omega,k_s>0$ and $x_R(R,q)$ is defined as in Theorem~\ref{theorem::kinematics}. If the potential function $\mathcal{U}\in\mathcal{P}_\mathcal{D}$ used in \eqref{q}-\eqref{D} and \eqref{control3} is {\it exp}-synergistic with synergistic gap exceeding $\delta$, then there exists a positive constant $\underline{k}_s$ such that for $k_s>\underline{k}_s>0$, the set
            $\bar{\mathcal{A}}_{s} :=\{(R,q,\omega,\tilde x_R)\in SO(3)\times\mathcal{Q}\times\mathbb{R}^3\times\mathbb{R}^3\mid\; R=I, \omega=0, \tilde x_R=0\}$
            is globally exponentially stable, where $\tilde x_R=x_{R,s}-x_R(R,q)$.
            \end{thm}
            \begin{remark}
            The control scheme in Theorem~\ref{theorem::dynamics::smooth} shows that the hybrid controller \eqref{control2} can be made \textit{continuous} without sacrificing global exponential stability. This is due to the integral action on the hybrid term $x_R(R,q)$ which moves the discontinuity one integrator away from the control input. The same remark can be noted for the input \eqref{du1}.
            \end{remark}
            The results in Theorem~\ref{theorem::kinematics} and Theorem~\ref{theorem::dynamics}~(Theorem~\ref{theorem::dynamics::smooth}) provide a unified method for the design of hybrid (continuous) attitude control algorithms for systems \eqref{dR:1} and \eqref{dR:2}, respectively, guaranteeing global exponential stability. It is worth pointing out that, in addition to the structurally simple expressions of the controllers proposed so far, the results in this subsection reduce the stabilisation problem of systems \eqref{dR:1} and \eqref{dR:2}, to the problem of finding appropriate {\it exp}-synergistic potential functions, in the sense of Definition~\ref{definition::expsynergism}, with some synergistic gap that can be specified using the control parameters.

\vspace{-0.05 in}
\subsection{Construction of Exp-Synergistic Potential Functions}% on $SO(3)\times\mathcal{Q}$}
\vspace{-0.1 in}
            In this work, we consider the following function %on $SO(3)$
            \begin{eqnarray}\label{V}
            V(R)=1-\sqrt{1-|R|_I^2}.
            \end{eqnarray}
            Note that $V$ in \eqref{V} is not differentiable on the set $\bar D_V=\{R\in SO(3)\mid\;|R|_I=1\}\subset SO(3)$ and is quadratic with respect to its gradient, see \citep{lee2012}. This makes it a good candidate for the construction of {\it exp}-synergistic potential functions on $SO(3)$ as shown in the following Proposition~\ref{proposition::V} proved in Appendix~\ref{proof::proposition::V}.
\vspace{-0.1 in}
            \begin{proposition}\label{proposition::V}
            Let $\mathcal{Q}=\{1,2,...,6\}$ be a finite index set. Consider the transformation $\Gamma: SO(3)\times\mathcal{Q}\to SO(3)$ defined as\vspace{-0.1 in}
            \begin{eqnarray}\label{Transformation2}
            \begin{array}{c}
            \Gamma(R,q)=R\mathcal{R}(2\arcsin\left(k |R|_I^2\right),u_q),
         \vspace{-0.1 in}
            \end{array}
            \end{eqnarray}
            where $0<k<1/\sqrt{2}$, $u_{m+3}=-u_m$, $m\in\{1,2,3\}$, and $\{u_1,u_2,u_3\}$ is an orthonormal set of vectors.
            Then, $V\circ\Gamma\in\mathcal{P}_\mathcal{D}$, where $V$ and $\Gamma$ are defined, respectively, in \eqref{V} and \eqref{Transformation2}, and $\mathcal{D}=\{(R,q)\in SO(3)\times\mathcal{Q}\mid\;|\Gamma(R,q)|_I\neq 1\}$. Moreover, $V\circ\Gamma$ is {\it exp}-synergistic with a gap $\delta > \bar{\delta}:=\frac{[-1+\sqrt{1+4k^2}]^{\frac{3}{2}}}{2\sqrt{6}k^2}$.
            \end{proposition}
\vspace{-0.1 in}
            \begin{remark}
            In Proposition \ref{proposition::V}, the transformation $\Gamma(R,q)$ can be seen as a perturbation to $R$ about the unit vector $u_q$ by an angle $2\arcsin\left(k |R|_I^2\right)$. This allows to \textit{stretch}/{\it compress} the manifold $SO(3)$ in order to move the singular points of $V(R)$ to different locations on $SO(3)$. Note that the transformation $\Gamma(R,q)$ in Proposition~\ref{proposition::V} is different from the one proposed in \cite{berkane2015construction} in the sense that \eqref{Transformation2} is specifically designed, taking into account \eqref{V}, to ensure the requirements of Definition~\ref{definition::synergism}. Consequently, the choice of the set $\mathcal{Q}$ in Proposition~\ref{proposition::V} leads to hybrid controllers \eqref{q}-\eqref{control3} with six possible modes. This choice is due to the nature of the singular manifold $\bar D_V\equiv \mathcal{R}(\pi,\mathbb{S}^2)$ which corresponds to the set of all rotations of angle $\pi$ around any vector in $\mathbb{S}^2$. Roughly speaking, to avoid this singular manifold one needs to apply angular warping in at least three different directions, where in each direction the warping angle could be positive or negative. This will be clear following the proof of the proposition.
            \end{remark}
            Now, consider the potential function $\mathcal{U}(R,q):= V\circ\Gamma(R,q)$ where $V(R)$ and $\Gamma(R,q)$ are defined in \eqref{V}-\eqref{Transformation2}. Noting that $\nabla\mathcal{U}(R,q)=\frac{\nabla|\Gamma(R,q)|_I^2}{2\sqrt{1-|\Gamma(R,q)|_I^2}}$, and using \eqref{psi}, \eqref{V}-\eqref{Transformation2}, and the result of \cite[Lemma 1 \& 3]{berkane2015construction}, the expression of the proportional-like term $x_R(R,q)$, used in \eqref{control1}-\eqref{control3}, can be derived as
            \begin{eqnarray}\label{beta2:1}
                x_R(R,q) =\frac{1}{8}\Theta_q^{\top}(R)\frac{\psi(\Gamma(R,q))}{\sqrt{1-|\Gamma(R,q)|_I^2}},
            \end{eqnarray}
           \vspace{-0.1 in}
           with
           \begin{eqnarray}\label{beta2:2}
           \Theta_q(R)=\mathcal{R}^\top(2\arcsin\left(k |R|_I^2\right),u_q)+\frac{ku_q\psi(R)^\top}{\sqrt{1-k^2 |R|_I^4}}.
            \end{eqnarray}
            \vspace{-0.1 in}
            Therefore, the following result holds.
            \begin{thm}\label{theorem:control:complete}
             Consider system \eqref{dR:1} (respectively \eqref{dR:2}) with the control input \eqref{control1} or \eqref{du1} (respectively \eqref{control2}) where $q\in\mathcal{Q}$ is generated by \eqref{q}-\eqref{D}, with $\mathcal{Q}$ and the potential function $\mathcal{U}(R,q):= V\circ\Gamma(R,q)$ given in Proposition~\ref{proposition::V}, and the term $x_R(R,q)$ is given by \eqref{beta2:1}-\eqref{beta2:2}. Pick the scalar $\delta$ in \eqref{q}-\eqref{D} such that $\delta>\bar{\delta}$, with $\bar{\delta}$ given in Proposition~\ref{proposition::V}.
%$\delta < \frac{[-1+\sqrt{1+4k^2}]^{\frac{3}{2}}}{2\sqrt{6}k^2}$ where $k$ is given in Proposition~\ref{proposition::V}.
             Then, the set $\mathcal{A}$ or $\mathcal{A}_{s}$ (respectively $\bar{\mathcal{A}}$) is globally exponentially stable. In addition, in the case of system \eqref{dR:2} with the control input \eqref{control3}, there exists a gain $k_s>\underline{k}_s>0$ such that $\bar{\mathcal{A}}_{s}$ is globally exponentially stable.
            \end{thm}
             The proof of Theorem~\ref{theorem:control:complete} follows by showing, from Proposition~\ref{proposition::V}, that the family of potential functions $\mathcal{U}(R,q):= V\circ\Gamma(R,q)$ satisfy all the conditions of Theorems~\ref{theorem::kinematics}-\ref{theorem::dynamics::smooth}. The results also show that the term $x_R(R,q)$ in \eqref{control1}, \eqref{control2} and \eqref{control3} is well defined in the flow set. In fact, since $k^2|R|_I^4\leq k^2<1/2$, the quantity in \eqref{beta2:2} is well defined. Also, using the fact that $|\Gamma(R,q)|_I^2\neq 1$, for all $(R,q)\in\mathcal{D}$, it follows that \eqref{beta2:1} is well defined on $C$ since $C\subseteq \mathcal{D}$ by virtue of the definition of the set $C$ and the fact that $\mathcal{U}(R,q)$ is {\it exp}-synergistic.
\vspace{-0.1 in}

\section{Application to Attitude Tracking on $SO(3)$}\label{Sec:appliation:examples}
\vspace{-0.1 in}
        In this section, we use the proposed approach to the attitude tracking control problem on $SO(3)$. Consider the attitude dynamics of a rigid body
            \begin{eqnarray}\label{dynamics}
            \dot{R}_\mathcal{B}^\mathcal{I}=R_\mathcal{B}^\mathcal{I}[\Omega]_\times,\qquad
            J\dot{\Omega}=[J\Omega]\times\Omega+\tau,
            \end{eqnarray}
            where $R_\mathcal{B}^\mathcal{I}\in SO(3)$ represents the orientation of the rigid body relative to $\mathcal{I}$, $\Omega\in\mathbb{R}^3$ is the angular velocity of the rigid body, $\tau\in\mathbb{R}^3$ is the input torque, and $J\in\mathbb{R}^{3\times 3}$ is the positive-definite inertia matrix of the rigid body with respect to $\mathcal{B}$. The control objective is to design a feedback control law $\tau$ such that $R_\mathcal{B}^\mathcal{I}\to R_d$ exponentially for arbitrary initial conditions, where $R_d(t)\in SO(3)$ is a desired attitude satisfying
            $\dot R_d=R_d[\Omega_d]_\times$, for some $R_d(0)\in SO(3)$ and a desired angular velocity $\Omega_d(t)\in\mathbb{R}^3$.

            Let $\tilde R_c := R_\mathcal{B}^\mathcal{I}R_d^{\top}$ and $\tilde\omega_c := R_d(\Omega-\Omega_d)$ denote, respectively, the attitude tracking error and the angular velocity tracking error. Also, consider the torque input in \eqref{dynamics}
            \begin{eqnarray}\label{torque}
            \tau=-[J\Omega]_\times\Omega+ J \big(\dot\Omega_d - [\Omega_d]_\times\Omega + R_d^{\top}u_c\big),
            \end{eqnarray}
             where $u_c$ is an additional input to be designed. Using \eqref{dynamics}, with \eqref{torque} and the property  $R^{\top}[x]_\times R=[R^{\top} x]_\times$ for every $x\in\mathbb{R}^3$ and $R\in SO(3)$, one can verify that
                 \begin{eqnarray}\label{kinematics_error}
                \dot{\tilde R}_c=\tilde R_c[\tilde\omega_c]_\times,\qquad \dot{\tilde\omega}_c=u_c,
                \end{eqnarray}
             which is similar to \eqref{dR:2}. Therefore, the following result follows immediately from Theorem~\ref{theorem:control:complete}.
\vspace{-0.1 in}
             \begin{corollary}\label{cor:traking1}
             Consider system \eqref{dynamics} with \eqref{torque}. Let $q_c\in\mathcal{Q}$ be generated by \eqref{q}-\eqref{D} with $\delta > \bar{\delta}$ and the potential function $\mathcal{U}(\tilde R_c,q_c):= V\circ\Gamma(\tilde R_c,q_c)$, where $\mathcal{Q}$,  $V(\tilde R_c)$, $\Gamma(\tilde R_c,q_c)$, and $\bar{\delta}$ are defined as in Proposition~\ref{proposition::V}. Also, let $u_c$ in \eqref{torque} be given as $u_c=-k_c x_R(\tilde R_c,q_c)-k_\omega\tilde\omega_c$, where $k_c,~k_\omega>0$ and $x_R(\tilde R_c,q_c)$ is given by \eqref{beta2:1}-\eqref{beta2:2}. Then, the set $\mathcal{A}_c:= \{(\tilde{R}_c,q_c, \tilde{\omega}_c\in SO(3)\times\mathcal{Q}\times \mathbb{R}^3\mid\;\tilde{R}_c=I, \tilde{\omega}_c=0\;\}$ is globally exponentially stable. If, instead, the continuous control input
             \begin{align*}
             u_c&=-k_cx_{R,s}-k_\omega\tilde\omega_c,\\
             \dot x_{R,s}&=-k_s(x_{R,s}-x_R(\tilde R_c,q_c)),
             \end{align*}
             is used, then there exists a gain $k_s>\underline{k}_s>0$ such that the set $\mathcal{A}_{c,s}:= \{(\tilde{R}_c,q_c, \tilde{\omega}_c, \tilde{x}_R)\in SO(3)\times\mathcal{Q}\times \mathbb{R}^3\times \mathbb{R}^3\mid\;\tilde{R}_c=I, \tilde{\omega}_c=0,\tilde{x}_{R}=0\;\}$, with $\tilde{x}_{R} = x_{R,s} - x_R(\tilde R_c,q_c)$, is globally exponentially stable.
             \end{corollary}
\vspace{-0.05 in}
\section{Conclusion}\label{conclusion}
\vspace{-0.1 in}
        We designed globally exponentially stabilizing hybrid feedback schemes for the attitude dynamics. Our approach relies on a central family of potential functions on $SO(3)$ that enjoy (as well as their gradients) the following properties: 1) being quadratic with respect to the Euclidean attitude distance, and 2) being synergistic with respect to the gradient's singular and/or critical points. An explicit construction of such potential functions achieving our control objectives, as well as a lower bound on the synergistic gap, are provided. In addition, we showed that it is possible to design continuous-time control schemes for the same problems without jeopardizing the global exponential stability results. The proposed hybrid scheme is applied to the attitude tracking problem on $SO(3)$ leading to global exponential stability; a result that was considered as the \textit{holy grail} for many years.

       % We designed globally exponentially stabilizing hybrid feedback schemes for the first and second order attitude dynamics. In addition, we showed that it is possible to remove jumps from the control input by adding a simple first order linear filter without jeopardizing the global exponential stability results.\\
%        The newly proposed hybrid feedback design relies on a new concept of {\it exp}-synergistic potential functions on $SO(3)$. An exp-synergistic potential function enjoys the property of  being (and its gradient) quadratic with respect to the Euclidean attitude distance and being synergistic with respect to the gradient possible singular points. An explicit construction of such a potential function as well as a lower bound on the synergistic gap are provided.\\
%        Furthermore, the proposed hybrid scheme is applied to the attitude tracking problem on SO(3) leading to global exponential stability; a result that was considered as the \textit{holy grail} for many years.
%\vspace{-0.1 in}
%%%%%%%%%%%%%%%%%%%%%%%%%%%%%%%55
%\vspace{-0.25cm}
\appendix
\section{Proof of Theorem \ref{theorem::kinematics}}\label{proof::theorem::kinematics}
\vspace{-.1 in}
            Let $\mathcal{U}\in\mathcal{P}_\mathcal{D}$, where $\mathcal{D}\subseteq SO(3)\times\mathcal{Q}$, be an {\it exp}-synergistic potential function with gap exceeding $\delta$. The closed loop system \eqref{dR:1} with \eqref{control1} along the flow set $C$ is written as
            \begin{eqnarray}\label{autonomuous::kinematics}
            \dot{R}=-k_cR\left[\psi(R^{\top}\nabla\mathcal{U}(R,q)))\right]_\times=-k_c\nabla\mathcal{U}(R,q).
            \end{eqnarray}
            Then, for all $(R,q)\in C$, the time derivative of $\mathcal{U}$ along \eqref{autonomuous::kinematics} is obtained as
            $
            \dot{\mathcal{U}}(R,q)=\langle\nabla\mathcal{U},\dot R\rangle_R\leq -k_c\alpha_3\;\mathcal{U}(R,q),
            $
            where we used \eqref{metric} and \eqref{condition::expsynergistic2} with \eqref{autonomuous::kinematics}. Moreover, for all $(R,q)\in D$, one has
            \begin{align*}
            \mathcal{U}(R,q^+)=\underset{m\in\mathcal{Q}}{\mathrm{min}}\;\mathcal{U}(R,m)&\leq\mathcal{U}(R,q)-\delta\leq e^{-\lambda_D}\;\mathcal{U}(R,q),
            \end{align*}
            where $\lambda_D=-\ln(1-\delta/\bar{\mathcal{U}})$ and $\bar{\mathcal{U}}:=\mathcal{U}(R(0,0),q(0,0))$ is the initial value of $\mathcal{U}$ which, since $\mathcal{U}$ is decreasing, represents the maximum of $\mathcal{U}(R,q)$ along the trajectories of \eqref{autonomuous::kinematics}. Consequently, we can apply the comparison principle \citep[Lemma C.1]{cai2009characterizations} to conclude that, for all $(t,j)\in\mathrm{dom}\;(R,q)$,$$
            \mathcal{U}(R(t,j),q(t,j))\leq\exp\left(-\lambda(t+j)\right)\mathcal{U}(R(0,0),q(0,0)),$$
            where $\lambda=\min\{k_c\alpha_3,\lambda_D\}>0$. Finally, since in view of \eqref{condition::expsynergistic1} $\;\mathcal{U}$ is quadratic with respect to $|R|_I$ it follows that the set $\mathcal{A}$ is globally exponentially stable.
\section{Proof of Theorem \ref{theorem::dynamics}}\label{proof::theorem::dynamics}
%\vspace{-.3cm}
Let $\mathcal{U}\in\mathcal{P}_\mathcal{D}$ be an {\it exp}-synergistic potential function with gap exceeding $\delta$ where $\mathcal{D}\subseteq SO(3)\times\mathcal{Q}$. The closed loop system \eqref{dR:2} with \eqref{control2} along $C$ is given by
            \begin{eqnarray}\label{autonomuous::dynamics1}
            \dot{R}=R[\omega]_\times,\qquad
            \dot{\omega}=-k_\omega\omega-k_cx_R(R,q),
            \end{eqnarray}
            which is an autonomous system.
           The angular velocity $\omega$ obeys dynamics of a first order linear filter with $x_R(R,q)$ as an input. It follows that
\begin{align*}
\|\omega(t,j)\|&\leq\|\omega_0\|e^{-k_\omega t}+\frac{k_c}{k_\omega}\underset{(t,j)\in\mathbb{R}^+\times\mathbb{N}}{\mathrm{sup}}\|x_R(R,q)\|\\
				   &\leq B_\omega:=\|\omega_0\|+\frac{k_c}{k_\omega}\sqrt{\alpha_4/2},
\end{align*}
where we have used the fact that, in view of \eqref{condition::expsynergistic2}, one has $\|x_R(R,q)\|^2\leq\alpha_4/2$. To show exponential stability, we consider the following Lyapunov function candidate
            \begin{eqnarray}\label{W}
            \mathcal{W}(R,q,\omega)=\frac{k_c}{2}\mathcal{U}(R,q)+\frac{1}{2}\|\omega\|^2+c\;\omega^{\top}x_R(R,q),
            \end{eqnarray}
            for some positive constant $c$. In view of \eqref{condition::expsynergistic1}-\eqref{condition::expsynergistic2}, one has
            $
            \frac{1}{2}z_1^{\top}M_1z_1\leq\mathcal{W}\leq\frac{1}{2}z_1^{\top}M_2z_1,
            $
            where $z_1=\left[|R|_I\;\|\omega\|\right]^\top$, and
            \begin{equation*}
            M_1=\begin{bmatrix}
            k_c\alpha_1&-c\sqrt{\frac{\alpha_4}{2}}\\
            -c\sqrt{\frac{\alpha_4}{2}}&1
            \end{bmatrix},\qquad
            M_2=\begin{bmatrix}
            k_c\alpha_2&c\sqrt{\frac{\alpha_4}{2}}\\
            c\sqrt{\frac{\alpha_4}{2}}&1
            \end{bmatrix}.
            \end{equation*}
            For all $(R,q)\in C$, the time derivative of $\mathcal{W}$ along \eqref{autonomuous::dynamics1} is given by\vspace{-0.7cm}
            {\small
            \begin{align*}\nonumber
            \dot{\mathcal{W}}=&k_cx_R(R,q)^\top\omega+\omega^\top(-k_\omega\omega-k_cx_R(R,q))\\
            &+c\omega^{\top}\dot{x}_R+cx_R^{\top}(-k_\omega\omega-k_cx_R)\\
            \leq&-ck_c\|x_R\|^2-k_\omega\|\omega\|^2+cD_x\|\omega\|^2+ck_\omega\|x_R\|\|\omega\|,
            \end{align*}}
            where we used the fact that $\dot x_R(R,q) = X(R,q)\omega$ for some $X(R,q)\in\mathbb{R}^{3\times 3}$. Also, since $SO(3)$ is a compact manifold, we know that there exists $D_x>0$ such that $\|X(R,q)\|_F\leq D_x$.  As a result, one can deduce that
             \begin{equation*}
             \dot{\mathcal{W}}\leq -\frac{1}{2}z_1^{\top}\begin{bmatrix}
            ck_c\alpha_3&-ck_\omega\sqrt{\frac{\alpha_4}{2}}\\
            -ck_\omega\sqrt{\frac{\alpha_4}{2}}&~~2(k_\omega-cD_x)\\
            \end{bmatrix}z_1= -\frac{1}{2}z_1^{\top}M_3z_1.
             \end{equation*}
             It can be verified that matrices $M_1,M_2$ and $M_3$ are positive definite for
                \begin{eqnarray*}
               0<c<\mathrm{min}\left\{2\sqrt{k_c/\alpha_4},\frac{4k_ck_\omega\alpha_3}{4D_xk_c\alpha_3+k_\omega^2\alpha_4}\right\},
                \end{eqnarray*}
             and therefore, one has
            $
            \dot{\mathcal{W}}\leq-\lambda_C\mathcal{W},$ for all $(R,q)\in C$, where $\lambda_C=\frac{\lambda_{\mathrm{min}}(M_3)}{\lambda_{\mathrm{max}}(M_2)}$. Moreover, for all $(R,q)\in D$, if the constant $c$ is chosen sufficiently small such that $c<\frac{k_c\delta}{8B_\omega\sqrt{\alpha_4/2}}$, one obtains
            \vspace{-0.8cm}{\small
            \begin{align*}
            \mathcal{W}(R,q^+,\omega)-\mathcal{W}(R,q,\omega)=&\frac{k_c}{2}\big(\mathcal{U}(R,g)-\mathcal{U}(R,q)\big)\nonumber\\
            & + c\omega^{\top}\big(x_R(R,g) - x_R(R,q)\big)\nonumber\\
            \leq & -\frac{k_c}{2}\delta+2c B_\omega\sqrt{\alpha_4/2} <-\frac{k_c}{4}\delta.
            \end{align*}}
Hence $\mathcal{W}$ is strictly decreasing over the jump set $D$. Also, since all signals are bounded, it follows that there exists $\bar{\mathcal{W}}$ such that $
            \mathcal{W}(R,q,\omega)\leq\bar{\mathcal{W}}
            $. Note also that for all $(R,q)\in D$, one has $\mathcal{U}(R,q)\geq\delta$ and therefore $\mathcal{W}(R,q,\omega,\tilde x_R)\geq k_c\delta/2$. This implies that $\mathcal{W}(R,g,\omega)\leq \exp(-\lambda_D)\mathcal{W}(R,q,\omega)$ where $\lambda_D=-\ln\left(1-\frac{k_c\delta}{4\bar{\mathcal{W}}}\right)$. Consequently, one can apply the comparison principle [Lemma C.1, \cite{cai2009characterizations}] to conclude that, for all $(t,j)\in\mathrm{dom}\;X$,
$$
\mathcal{W}(X(t,j))\leq\exp(-\lambda(t+j))\mathcal{W}(X(0,0)),
$$
where $\lambda=\mathrm{min}\{\lambda_C,\lambda_D\}$. This shows, according to \cite{teel2013lyapunov}, that the equilibrium $z=0$ is globally exponentially stable.
 \vspace{-0.4cm}
\section{Proof of Theorem \ref{theorem::dynamics::smooth}}\label{proof::theorem::dynamics::smooth}
%\vspace{-.3cm}
Let $\mathcal{U}\in\mathcal{P}_\mathcal{D}$ be an {\it exp}-synergistic potential function with gap exceeding $\delta$ and where $\mathcal{D}\subseteq SO(3)\times\mathcal{Q}$. The closed loop system \eqref{dR:2} with \eqref{control2} along $C$ is given by
            \begin{align}\nonumber
            \dot{R}&=R[\omega]_\times,\\\nonumber
            \dot{\omega}&=-k_cx_{R,s}-k_\omega\omega,\\
            \dot x_{R,s}&=-k_s(x_{R,s}-x_R(R,q)),\label{autonomuous::dynamics2}
            \end{align}
            which is an autonomous system. The continuous variable $x_{R,s}$ is the output of a first order linear system with $x_R(R,q)$ being a bounded discontinuous input. Therefore, $x_{R,s}$ is bounded and so is $\tilde x_R=x_{R,s}-x_R(R,q)$. Let us denote $B_{\tilde x}$ such a bound. Also using a similar argument one can show that the angular velocity $\omega$ is bounded by some $B_\omega$.
 Consider the following Lyapunov function candidate
\begin{align*}
\mathcal{W}_s=\frac{k_c}{2}\mathcal{U}(R,q)+\frac{1}{2}\|\omega\|^2+c_1x_R^\top\omega+\frac{c_2}{2}\|\tilde x_R\|^2,
\end{align*}
where $c_1$ and $c_2$ are some positive constants. It is clear that, if $c_1<2\sqrt{k_c/\alpha_4}$, one has $\mathcal{W}_s$ is positive definite with respect to $\mathcal{A}_{2,s}$. For all $(R,q)\in C$, the time derivative of $\mathcal{W}_s$ along the trajectories \eqref{autonomuous::dynamics2} is given by
{\small
\begin{align*}
            \dot{\mathcal{W}}_s=&k_cx_R(R,q)^\top\omega+\omega^\top(-k_\omega\omega-k_cx_{R,s})+c_1\omega^\top X(R,q)\omega\\
            &+c_1x_R^\top(-k_\omega\omega-k_cx_{R,s})+c_2\tilde x_R^\top (-k_s\tilde x_R-X(R,q)\omega)\\
            \leq&-c_1k_c\|x_R\|^2-k_\omega\|\omega\|^2-c_2k_s\|\tilde x_R\|^2+c_1k_c\|x_R\|\|\tilde x_R\|\\
            &+(k_c+c_2D_x)\|\omega\|\|\tilde x_R\|+ c_1k_\omega\|x_R\|\|\omega\|+c_1D_x\|\omega\|^2\\
            \leq&-\frac{1}{2}z_1^\top P_1z_1-\frac{1}{2}z_2^\top P_2z_2-\frac{1}{2}z_3^\top P_3z_3,
 \end{align*}}
 where
%{\small
 \begin{align*}
P_1=\begin{bmatrix}
 c_1k_c\alpha_3/2&-c_1k_\omega\sqrt{\alpha_4/2}\\
 -c_1k_\omega\sqrt{\alpha_4/2}&k_\omega-2c_1D_x
 \end{bmatrix},\quad\\
%%%%%%%%%%%%%%%%%
 P_2=\begin{bmatrix}
 c_1k_c\alpha_3/2&-c_1k_c\sqrt{\alpha_4/2}\\
 -c_1k_c\sqrt{\alpha_4/2}&c_2k_s
 \end{bmatrix},\quad
 %%%%%%%%%%%%%%%%%
\\ P_3=\begin{bmatrix}
k_\omega&-(k_c+c_2D_x)\\
-(k_c+c_2D_x)&c_2k_s
 \end{bmatrix},
 \end{align*}%}
and $z_1=[|R|_I,\;\|\omega\|]^\top$, $z_2=[|R|_I,\;\|\tilde x_R\|]^\top$ and $z_3=[\|\omega\|,\;\|\tilde x_R\|]^\top$. The matrices $P_1$ and $P_2$ are positive definite provided that
$$
0<c_1<\mathrm{min}\left\{\frac{k_ck_\omega\alpha_3}{k_\omega^2\alpha_4+2k_cD_x\alpha_3},\;\frac{c_2k_s\alpha_3}{k_c\alpha_4}\right\}.
$$
Moreover, $P_3$ is positive if one chose $
k_s>\underline{k}_s:=(k_c+c_2D_x)^2/c_2k_\omega.
$
In this case there exists $\lambda_C$ such that $\dot{\mathcal{W}}_s\leq-\lambda_C\mathcal{W}_s$ for all $(R,q)\in C$.
 Between jumps, we have for all $(R,q)\in D$
%{\small
\begin{align*}
&\mathcal{W}_s(R,g,\omega,\tilde x_R)-\mathcal{W}_s(R,q,\omega,\tilde x_R)\\
&\;\;\;\leq-\frac{k_c\delta}{2}+c_1x_R^\top\omega\left|_{(R,q)}^{(R,g)}\right.+\frac{c_2}{2}\|\tilde x_R\|^2\left|_{(R,q)}^{(R,g)}\right.\leq\frac{k_c\delta}{4},
\end{align*}
where we further imposed that $c_1<k_c\delta/16\sqrt{\alpha_4/2}B_\omega$ and $c_2<k_c\delta/8B_{\tilde x}^2$. Hence, $\mathcal{W}_s$ is strictly decreasing over jumps. Since all signals are bounded, there exists $\bar{\mathcal{W}}$ such that $\mathcal{W}_s(R,q,\omega,\tilde x_R)\leq\bar{\mathcal{W}}$. Note also that for all $(R,q)\in D$, one has $\mathcal{U}(R,q)\geq\delta$ and therefore $\mathcal{W}_s(R,q,\omega,\tilde x_R)\geq k_c\delta/2$. Hence one can write $\mathcal{W}_s(R,g,\omega,\tilde x_R)\leq\exp(-\lambda_D)\mathcal{W}_s(R,q,\omega,\tilde x_R)$ where $\lambda_D=-\mathrm{ln}\left(1-k_c\delta/4\bar{\mathcal{W}}\right)$. Consequently, it follows that $\mathcal{A}_{2,s}$ is globally exponentially stable.
%\begin{remark}
%To compute the value of the lower bound $\underline{k}_s$ one needs to have access to $D_x$ which is the upper bound of $X(R,q)$. To avoid computing $D_x$ we may use the control
%\begin{align*}
%\nonumber
%u&=-k_cx_{R,s}-k_\omega\omega,\\
%\dot x_{R,s}&=X(R,q)\omega-k_s(x_{R,s}-x_{R}(R,q)),\label{x_smooth3.1}
%\end{align*}
%which compensate for the non-linearity $X(R,q)\omega$. Therefore we may end up with the condition
%\begin{align}
%k_s&>\underline{k}_s:=\frac{k_c^2}{c_2k_\omega}\\
%c_2&=\frac{k_c\delta}{8B_{\tilde x}^2},
%\end{align}
%which is available for computation. The user is given the choice between implementing $X(R,q)\omega$ or computing the upper bound $D_x$ of the matrix $X(R,q)$.
%\end{remark}
\section{Proof of Proposition \ref{proposition::V}}\label{proof::proposition::V}
Let $(\eta,\epsilon)$ and $(\eta_q,\epsilon_q)$ be the quaternion representations of the attitude matrices $R$ and $\Gamma(R,q)$, respectively. The unit quaternion associated to $\mathcal{R}(2\arcsin(k|R|_I^2),u_q)$  is given by $(\sqrt{1-k^2|R|_I^4},k|R|_I^2u_q)$. Note that, it can be checked that $|R|_I^2=\|\epsilon\|^2$. Therefore, in view of \eqref{Transformation2} and using the quaternion multiplication rule, one obtains
\begin{align}\label{etaq}
\eta_q&=\eta\sqrt{1-k^2\|\epsilon\|^4}-k\|\epsilon\|^2\epsilon^{\top}u_q,\\
\epsilon_q&= k\eta\|\epsilon\|^2u_q+\sqrt{1-k^2\|\epsilon\|^4}\epsilon+k\|\epsilon\|^2\epsilon\times u_q.
\end{align}
This leads to write
\begin{align*}
\|\epsilon_q\|^2&=\|\epsilon\|^2+k^2\|\epsilon\|^4\eta^2-\cos^2(\varphi_q)k^2\|\epsilon\|^6+\\
						&\hspace{0.3cm}2k\eta\|\epsilon\|^3\sqrt{1-k^2\|\epsilon\|^4}\cos(\varphi_q)
\end{align*}
where $\varphi_q$ is the angle between $\epsilon$ and $u_q$. Using the fact that $|\eta|\cdot\|\epsilon\|\leq\frac{1}{2}$, it follows that
\begin{align*}
\|\epsilon_q\|^2&\leq\|\epsilon\|^2[1+k+k^2/4].
\end{align*}
Moreover, since $k<\frac{1}{\sqrt{2}}$ it is possible to show that\\
\begin{align*}
\|\epsilon_q\|^2&\geq\|\epsilon\|^2\sqrt{1-k^2\|\epsilon\|^4}\left[\sqrt{1-k^2\|\epsilon\|^4}-2k|\eta|\cdot\|\epsilon\|\right]\\
	 &\geq \|\epsilon\|^2\sqrt{1-k^2\|\epsilon\|^4}\left[\sqrt{1-k^2}-k\right]\\
	 &\geq \|\epsilon\|^2\left[1-k^2-k\sqrt{1-k^2}\right].
\end{align*}
On the other hand, one has $|\Gamma(R,q)|_I^2=\|\epsilon_q\|^2$ and
\begin{align*}
\frac{1}{2}|\Gamma(R,q)|_I^2\leq1-\sqrt{1-|\Gamma(R,q)|_I^2}&\leq|\Gamma(R,q)|_I^2,
\end{align*}
thanks to the fact that $|\Gamma(R,q)|_I\in[0,1]$. Subsequently, it follows that $\alpha_1|R|_I^2\leq V\circ\Gamma(R,q)\leq\alpha_2|R|_I^2$ with $\alpha_1=\left[1-k^2-k\sqrt{1-k^2}\right]/2$ and $\alpha_2=[1+k+k^2/4]$. It is easy to check that $\alpha_1$ and $\alpha_2$ are strictly positive for all $0\leq k<1/\sqrt{2}$.\\
Next, we show that the potential function $V\circ\Gamma$ satisfies \eqref{condition::expsynergistic2} of Definition \ref{definition::expsynergism}. By \cite[Lemma 3]{berkane2015construction}, one has $\det(\Theta_q(R))\neq 0$ for all $R\in SO(3)$. This implies that the matrix $\Theta_q(R)\Theta_q(R)^{\top}$ is full rank and positive definite for all $(R,q)\in SO(3)\times\mathcal{Q}$. Let us denote $\underline{\lambda}_\Theta^*, \bar\lambda_\Theta^*>0$ the smallest, respectively greatest, eigenvalue of $\Theta_q(R)\Theta_q(R)^{\top}$ for all $SO(3)\times\mathcal{Q}$. Then, in view of \eqref{beta2:1}-\eqref{beta2:2}, one has
\begin{equation*}
\begin{split}
\frac{\underline{\lambda}_\Theta^*}{64}\frac{\|\psi\circ\Gamma\|^2}{1-|\Gamma|_I^2}\leq\|\psi\left(R^{\top}\nabla(V\circ\Gamma)\right)
					  \|^2\leq \frac{\bar\lambda_\Theta^*}{64}\frac{\|\psi\circ\Gamma\|^2}{1-|\Gamma|_I^2}.
\end{split}
\end{equation*}
Furthermore, it is straightforward to verify that
$
\|\psi(X)\|^2=4\Phi_I(X)(1-\Phi_I(X)),\;\forall X\in SO(3).
$
It follows that
\begin{equation*}
\begin{split}
\frac{\underline{\lambda}_\Theta^*}{32}|\Gamma|_I^2\leq\|\psi\left(R^{\top}\nabla(V\circ\Gamma)\right)
					  \|^2\leq \frac{\bar\lambda_\Theta^*}{32}|\Gamma|_I^2.
\end{split}
\end{equation*}
Finally, one concludes that the potential function $V\circ\Gamma$ satisfies condition \eqref{condition::expsynergistic2} with $\alpha_3=\underline{\lambda}_\Theta^*\alpha_1/16$ and $\alpha_4=\bar\lambda_\Theta^*\alpha_2/32$.\vspace{-.3cm}

Let us show that the potential function $V\circ\Gamma$ satisfies \eqref{condition::expsynergistic3} in Definition~\ref{definition::expsynergism}. The potential function $V\circ\Gamma$ is not differentiable only on the set
$
\bar{\mathcal{D}}:=\left\{(R,q)\in SO(3)\times\mathcal{Q}\mid\;|\Gamma(R,q)|_I=1\right\}.
$
Let us evaluate the transformations $\Gamma(R,q)$ and $\Gamma(R,m)$, defined in \eqref{Transformation2}, for some $q,m\in\mathcal{Q}$ on the set $\bar{\mathcal{D}}$. Let $Q=(\eta,\epsilon)$, $Q_{q}=(\eta_{q},\epsilon_{q})$ and $Q_{m}=(\eta_{m},\epsilon_{m})$ be the quaternion representation of the attitude $R$, $\Gamma(R,q)$ and $\Gamma(R,m)$, respectively. Let us define the set
$
\bar{\mathcal{D}}_Q=\left\{((\eta,\epsilon),q)\in\mathbb{Q}\times\mathcal{Q}\mid\;\eta_q=0,\;\epsilon_q\in\mathbb{S}^2\right\}.
$
Note that the set $\bar{\mathcal{D}}_Q$ represents a double cover of the set $\bar{\mathcal{D}}$. It is not difficult to show that
$
|R|_I^2=1-\eta^2,
$
for all $R\in SO(3)$. Hence, one has
\begin{eqnarray}
\begin{array}{ll}
V\circ\Gamma(R,m)=1-|\eta_m|,\\
V\circ\Gamma(R,q)=1,\;\;\;\forall (R,q)\in\bar{\mathcal{D}}.
\end{array}
\end{eqnarray}
Therefore, if one guarantees that ${\mathrm{max}}_{m\in\mathcal{Q}}\;|\eta_m|>\delta$, for some $\delta>0$, then one has the condition $C\subseteq\mathcal{D}$ satisfied. In view of \eqref{etaq} and for all $(R,q)\in\bar{\mathcal{D}}$, one obtains $\eta_m-\eta_q=\eta_m=k\|\epsilon\|^2\epsilon^{\top}(u_q-u_m).
$ Moreover, since $u_{m+3}=-u_m, m\in\{1,2,3\}$ there exist always three indices $m_i\in\mathcal{Q}, i\in\{1,2,3\},$ such that $u_{m_i}$ is orthonormal to $u_{m_j},$ for $i\neq j$ and $\epsilon^\top u_{m_i}$ has opposite sign to $\epsilon^\top u_q$. Therefore, one has
\begin{align*}
\max_{m\in\mathcal{Q}}|\eta_m|&=k\|\epsilon\|^2\max_{m\in\mathcal{Q}}|\epsilon^\top(u_q-u_m)|,\\
&\geq k\|\epsilon\|^2\max_{m_i\in\{1,2,3\}}|\epsilon^\top u_{m_i}|\geq k\|\epsilon\|^3/\sqrt{3}.
\end{align*}
Furthermore, for all $(Q,q)\in\bar{\mathcal{D}}_Q$, equation \eqref{etaq} implies that
\begin{equation*}
\sqrt{1-\|\epsilon\|^2}\sqrt{1-k^2\|\epsilon\|^4}=|k|\|\epsilon\|^3|\cos(\varphi_q)|.
\end{equation*}
The above equation reads $f(\|\epsilon\|^2)=0$ where
$
f(x):=g(x)+\sin^2(\varphi_q)k^2x^3,
$
and $g(x):=1-x-k^2x^2$. It is easy to verify that $f(x)$ and $g(x)$ are decreasing on the interval $[0,1]$ for all $\varphi_q\in\mathbb{R}$. Therefore, since $f(x)\geq g(x)$, the solution $x_f\in[0,1]$ to equation $f(x_f)=0$ is greater than or equals to $x_g\in[0,1]$, with $g(x_g)=0$. Thus, it is clear that $f(\|\epsilon\|^2)=0$ implies that
$
\|\epsilon\|^2\geq(-1+\sqrt{1+4k^2})/2k^2.
$
Finally, one concludes that
 $\underset{m\in\mathcal{Q}}{\mathrm{max}}\;|\eta_m|\geq [-1+\sqrt{1+4k^2}]^{\frac{3}{2}}/2\sqrt{6}k^2.
$
%%%%%%%%%%%%%%%%%%%%%%%%%%%%%%%%%%%%%%%%55555
\bibliographystyle{apalike}
\bibliography{Hybrid}
\end{document}